\documentclass[10pt,a4paper,oneside]{article}
\pdfoutput=1
\usepackage[utf8]{inputenc}
\usepackage{amsmath}
\usepackage{amsfonts}
\usepackage{amssymb}

\usepackage{tikz-cd}
\usepackage[square,numbers]{natbib} 

\DeclareMathOperator{\Cpct}{Cap}
\DeclareMathOperator{\Co}{Co}
\def\eps{\varepsilon}
\def\epsi{\varepsilon} 
\def\Om{\Omega}
\newcommand{\BB}{{\mathfrak{B}}}
\newcommand{\sm}{\setminus}  
\newcommand{\C}{\mathbb{C}}
\newcommand{\Cstar}{{\mathbb{C^*}}}
\newcommand{\D}{\mathbb{D}}
\newcommand{\aphi}{\phi_{a,b}} 
\newcommand{\eul}{\textup{e}}
\newcommand{\vv}{\vert}
\newcommand{\phat}{\widehat{p}}
\newcommand{\gammahat}{\widehat{\gamma}}
\newcommand{\aphifn}{{\aphi \circ f^n}}
\newcommand{\aphifm}{{\aphi \circ f^m}}
\newcommand{\aphifk}{{\aphi \circ f^k}}

\newcommand{\qhat}{\widehat{q}}

\usepackage{amsthm} 

\newtheorem{theorem}{Theorem}[section] 
\newtheorem{definition}[theorem]{Definition}
\newtheorem{lemma}[theorem]{Lemma}
\newtheorem{proposition}[theorem]{Proposition} 
\newtheorem{corollary}[theorem]{Corollary}

\newtheorem{remark}[theorem]{Remark}

\title{Conformal Equivalence of  Measures and Dynamics of Orthogonal Polynomials\footnote{This is an Accepted Manuscript of an article published by Taylor \& Francis in 
\textbf{Journal of Difference Equations and Applications}, 25:8, 1063-1081, DOI: 10.1080/10236198.2019.1648446}}
\author{Signe Emalia Jensen\textsuperscript{1}   
and  Carsten Lunde Petersen\textsuperscript{2} \\ \small \textsuperscript{1} Northwestern University; \textsuperscript{2} Roskilde University \normalsize }
\date{ }

\begin{document}

\maketitle

\begin{abstract}

We introduce a notion of asymptotically orthonormal polynomials for a Borel measure $\mu$ with compact nonpolar support in $\C$. 
Such sequences of polynomials  have similar convergence properties of the sequences of Julia sets and filled Julia sets to those for sequences of orthonormal polynomials, see also \citep{lunde}.
We give examples of measures for which the  monic orthogonal polynomials are asymptotically orthonormal. 
Combining this with observations on conformal invariance of orthogonal polynomials we explore the measure dependency of the associated dynamics of  orthogonal polynomials. 
Concretely, we study the dynamics of sequences of asymptotically orthonormal polynomials for the pullback measure $\phi^\ast(\mu)$ under affine mappings $\phi$. 
We prove that the sequences of Julia sets and filled Julia sets of  affine deformations of sequences of asymptotically orthonormal polynomials for $\mu$ also have the same convergence properties as the Julia sets and filled Julia sets of the orthonormal polynomials.
This leads to theorems on the convergence properties of affine deformations of the family of iterates of any fixed monic  centered  polynomial and, in the case the polynomial is hyperbolic, on the corresponding family of affine parameter spaces. 

\end{abstract}

\section{Introduction}

Let $\mathfrak{B}$ denote the set of Borel probability measures on $\C$ with compact nonpolar support and let $\mu\in\mathfrak{B}$. 
We study sequences $\{q_n(\mu;z)\}_n$ of orthogonal polynomials for $\mu$.  
More precisely let $\{ p_n (\mu ; z) \}_n $ denote the unique orthonormal polynomial  sequence in $L^2 (\mu)$ with 
\begin{equation}
\label{eq:sequenceelement} 
p_n ( \mu ; z) = \gamma_n z^n + \textrm{lower order terms}
\end{equation}
and $\gamma_n>0$. Then for any sequence $\{a_n\}_n\subset\C$ of non-zero complex numbers the 
polynomial sequence $\{q_n(z) = q_n(\mu;z) := a_np_n(\mu;z)\}_n$ 
is an orthogonal sequence in $L^2 (\mu)$. 
Moreover as is standard in the literature we denote by $S(\mu)$ the compact support of $\mu$, 
by $\Om = \Om _ \mu $
the unbounded connected component of $S(\mu)$ and by $g_\Om$ the Green's function for $\Om$ with pole at infinity. 
In this paper we study the variation of such sequences as the measure $\mu$ is pulled back by an affine map and we are in particular interested in the dynamics of such polynomials. 

Polynomial dynamics is the study of dynamical systems given by iteration 
of a fixed polynomial $q$ of degree $n\geq 2$. That is, we consider orbits 
$z_0, z_1, \ldots , z_k, z_{k+1}, \ldots$ where 
$$
z_{k+1} = q(z_k) = q\circ q(z_{k-1}) = q^{\circ k}(z_0).
$$
We shall review  the basic notions of polynomial dynamics used here, 
the interested reader may consult \cite{CarlesonandGamelin} or \cite{milnor} for more comprehensive introductions to the theory of polynomial 
and more generally holomorphic dynamics. 

The dynamics of $q$ naturally divides $\C$ into two completely invariant complementary sets, the open Fatou set $F(q)$ with `tame' dynamics and the compact Julia set $J(q)=\C\sm F(q)$ with chaotic dynamics. Technically, the Fatou set is the set points $z$ 
possessing a neighbourhood on which the family of iterates $\{q^{\circ k}\}_k$ 
form a normal family in the sense of Montel, so that the Julia set is the set of points $z$ 
such that the family of iterates $\{q^{\circ k}\}_k$ is not normal on any neighbourhood. 
The Julia set is also the closure of the set of repelling periodic points for $q$, 
i.e. points $z$ such that $q^{\circ k}(z) = z$ for some $k\geq 1$ and $|(q^{\circ k})'(z)|> 1$.
The unbounded connected component $\Om_q$ of $F_q$ consists of points $z$ such 
that $q^{\circ k}(z) \to\infty$ as $k\to\infty$ and is called the basin of attraction for $\infty$. 
The complement $K(q)=\C\setminus\Om_q$ is called the filled Julia set of $q$. 
The common boundary of the two sets is precisely the Julia set of $q$. 
The Julia set is Dirichlet regular and the Green's function $g_q$ for $\Om_q$ extended 
by zero to $K(q)$ can be obtained as the locally uniform limit
$$
g_q(z) = \lim_{k\to\infty}\frac{1}{n^k}\log^+|q^{\circ k}(z)|
$$
where $\log^+x := \max\{0, \log x\}$ for any $x\geq 0$.

The dynamics of orthogonal polynomials is described in 
 the paper ``Julia Sets of Orthogonal Polynomials'' \citep{lunde}. 
In \citep{lunde} it is shown that the Julia sets and filled Julia sets for the sequence $\{ p_n (\mu ; z) \}$ of orthonormal polynomials is related to the support $S (\mu)$ of $\mu$.  
 The study relies on the monograph ``General Orthogonal Polynomials'' \citep{stahltotik}, where potential theoretic and measure theoretic properties of the support of a measure $\mu$ are related to the potential theoretic and measure theoretic properties
of the sequence of orthonormal polynomials defined by $\mu$.
The relations between the dynamics of $p_n (\mu ; z)$ and the support $S (\mu)$ are proved in \citep{lunde} using estimates on various Green's functions. 

In the present paper we study the dynamics of a broader class of orthogonal polynomials. We introduce here a notion of  asymptotically orthonormal polynomials.
 
\begin{definition}
\label{def:asymptoticallyorthonormal}
Let  $\mu\in\mathfrak{B}$ and let $\{a_n\}_n $ be a sequence with 
$a_n \in \C \sm \{ 0 \}$.
An asymptotically orthonormal sequence for $\mu$ is an orthogonal sequence 
$\{q_n (\mu ; z) = a_n p_n (\mu ; z)\}_n$ where $\lim_{n \rightarrow \infty} \sqrt[n]{ \vv a_n \vv} = 1$.
\end{definition}

The term asymptotically orthonormal derives from the fact, which we are to prove here, that 
such sequences behave in many aspects as the sequences of orthonormal polynomials.
If $\vv a_n \vv = 1$ for all $n$, then $\{q_n (\mu ; z)\}_n$ 
is an orthonormal sequence for $\mu$, although not normalized to have leading coefficient a positive real. We denote  such  orthonormal sequences $\{ \phat_n (\mu ; z) \}$. 
 
The natural notion of dynamical equivalence for polynomials is conformal or affine equivalence: 
two polynomials $p$ and $q$ are said to be conformally equivalent, if there is a conformal equivalence of $\C$ to $\C$, i.e.~an affine map $\phi(z) = \aphi(z) = az + b$, 
where $a \not= 0$ such that 
$$
p = \phi\circ q\circ \phi^{-1}.
$$
Note that necessarily the two polynomials $p$ and $q$ have the same degree. 
Moreover $\phi$ maps orbits of $q$ to orbits of $p$, $\Om_q$ to $\Om_p$, $K(q)$ to $K(p)$ and $J (q )$ to $J (p )$. 
Conformal equivalence of polynomials as dynamical systems is of great use in reducing 
the complexity of parameter spaces, while keeping all possible dynamics.

When studying the variation of orthogonal and orthonormal polynomials with the measure $\mu\in\BB$ we may similarly reduce the complexity by defining two measures 
$\mu,\nu\in\BB$ to be conformally equivalent, written $\mu\sim\nu$, 
if there is an affine mapping $\phi$ as above such that $\nu = \phi_*(\mu)$ 
is the conformal pushforward of $\mu$, 
i.e. $\nu(U) = \mu(\phi^{-1}(U))$ for any Borel set $U\subset\C$. 
The natural operation on measures is the pushforward by functions or maps, but since affine maps $\phi$ are biholomorphic we may as well consider 
the pushforward by $\phi^{-1}$, which technically is a pullback. 
In general the pullback of a measure by some map involves a choice. 
But in this case, where $\phi$ is a homeomorphism, it is simply the inverse operation. 
So in order to alleviate notation we use  the pullback 
by $\phi$ of measures instead of the pushforward by $\phi$, 
i.e. consider $\nu = \phi^{-1}_*(\mu) = \phi^\ast (\mu)$.

It immediately follows that the support $S (\nu)$ is the corresponding affine transformation 
of the support of $\mu$, i.e. $S(\mu) = \phi (S (\nu))$.

Moreover, for the conformal pullback measure $\nu$ of $\mu$  we have that 
for any orthogonal sequence $\{q_n(\mu ; z)\}_n$
\begin{equation}
\label{eq:polynomialtransformationq}
q_n (\nu ; z ) := q_n (  \mu ; \phi(z) ) 
\end{equation} 
is an orthogonal sequence for $\nu$. 
Equivalently  $ q_n \left(  \mu ; z \right) := q_n (\nu ; \phi^{-1}(z)  ) $ 
is an orthogonal sequence for $\mu$ whenever $\{q_n (\nu ; z )\}_n $ 
is an orthogonal sequence for $\nu$. 
And furthermore  
\begin{equation} \label{eq:polynomialtransformation}
\phat_n (\nu ; z ) := \phat_n (  \mu ; \phi  (z) ) 
\end{equation} 
is an orthonormal sequence for  $\nu$  whenever $\phat_n (\mu ; z ) $ is an orthonormal sequence for $\mu$. However we remark that $p_n (\nu , z ) \neq p_n (\mu , \aphi (z))$  unless $a$ is a positive real. 

The pullback of sequences of orthogonal and orthonormal polynomials have new dynamics,  
since $q_n(\mu ; z)$ and $q_n(\nu ; z) = q_n(\mu ;\phi(z))$ are not conformally conjugate. 
That is, being equivalent from the point of view of orthogonal polynomials is different from being equivalent from the point of view of dynamical systems. 
Thus the orthonormal sequence of polynomials for even a nearby 
conformally equivalent measure has new dynamics. 
We can view this as saying the results in \cite{lunde} have new content for 
each measure $\nu$ conformally equivalent to $\mu$. 

In this paper we study the dynamics of both variations of the orthonormal polynomials 
for $\mu\in\BB$,  the asymptotically orthonormal sequences of polynomials and the pullback.

We observe furthermore that $q_n(\nu ; z)$ is conformally conjugate to 
$\phi\circ q_n(\mu ; z)$, see Lemma~\ref{prop:conformalconjugate}, 
which leads us to study polynomial sequences of the form $\{\phi \circ p_n(\mu, z)\}_n$ 
and $\{\phi\circ q_n(\mu, z)\}_n$. 

The following Main Theorem uses the notions of minimal carrier Green's function $g_\mu$ as well as the notion of limits of sequences of compacts sets. These will be described in the following section.

\begin{theorem}
\label{thm:mainq}
Let $\mu \in \BB$, let $\phi$ be affine and let $\nu =  \phi^{-1} _\ast (\mu ) = \phi^\ast (\mu)$.
Suppose $q _n(z) = a_n p_n (\mu ; z)$ is an orthogonal polynomial sequence for $\mu$.
Then the support $S(\mu)$ and the Julia sets $J( {\phi\circ q_n} )$ 
and filled Julia sets $K ( {\phi\circ q_n} )$ have the same relative position as 
the support $S(\nu)$ of the measure $\nu$, the Julia sets $J({q_n\circ\phi})$ 
and filled Julia sets $K ( {q_n\circ\phi})$. 
Furthermore , if $\mu \in \BB_+$ and 
the sequence $\{q_n\}_n$ is  asymptotically orthonormal for $\mu$, then
\begin{equation}
\label{mainq1}
\overline{J _ \mu \sm E_{\mu}} \subseteq \liminf_{n\rightarrow\infty} J ({\phi\circ q_n} )\subseteq \limsup_{n\rightarrow\infty} K ({\phi\circ q_n}) \subseteq \textup{Co}(S(\mu)) ,
\end{equation} 
where  $J_ \mu$ is the outer boundary of the support $S(\mu)$, 
while $E_\mu$ denotes the exceptional set $\{ z \in S(\mu) \mid g_\mu (z) > 0 \}$.
Moreover let $V_\epsi := \{z\in\C\; |\; g_{\Om  _ \mu }(z) \geq \epsi\}$. 
Then 
\begin{equation}
\label{mainq2}
\forall\; \epsi>0\;:\qquad \lim_{n\to\infty} \Cpct(V_\epsi\cap K({\phi\circ q_n})) = 0 .
\end{equation}
\end{theorem}
Here $\mathfrak{B}_+$ denotes the set of measures 
$ \{ \mu \in \BB \mid \displaystyle{\liminf_{n\rightarrow\infty}} \Cpct (K({p_n (\mu ; z)})) > 0 \}$.

For a measure $\mu\in\mathfrak{B}$ the sequence of othonormal polynomials 
$\{p_n (\mu ; z)\}_n$ is said to have regular $n$th root asymptotic behavior, 
if and only if 
$$\lim_{n \rightarrow \infty} \left\vert p_n (\mu ; z) \right\vert^{\frac{1}{n}} 
= \eul^{g_\Omega (z)}
$$ locally uniformly in $\mathbb{C} \setminus \textup{Co}(S(\mu))$. 
Equivalently 
\begin{equation}
\label{eq:nthrootreg1}
\lim_{n \rightarrow \infty} \left(\frac{1}{\gamma_n}\right)^{\frac{1}{n}} = \Cpct(S(\mu)),
\end{equation}
where $\Cpct(S(\mu))$ denotes the logarithmic capacity of $S(\mu)$. 
We shall abuse notation and say that such a measure $\mu$ 
is $n$th root regular. Evidently such a measure belongs to $\mathfrak{B}_+$.

\begin{theorem}
\label{thm:mainregular}
Let $\mu$ be  $n$-th root regular, let $\phi$ be affine and  
let $q _n(z) = a_n p_n (\mu ; z)$ be an asymptotically orthonormal sequence for  $\mu$.
Then 
\begin{equation}
\label{mainq1nthroot}
\overline{J_\mu \sm E} \subseteq \liminf_{n\rightarrow\infty} J ({\phi\circ q_n}) \subseteq \limsup_{n\rightarrow\infty} K({\phi\circ q_n}) \subseteq \textup{Co}(S(\mu)) ,
\end{equation} 
where  $J_\mu$ is the outer boundary of the support $S(\mu)$, 
and $E$ denotes the $F_\sigma$ and polar exceptional set 
$\{ z \in S(\mu) \mid g_{\Omega_\mu} (z) > 0 \}$. In particular 
\begin{equation}
\label{mainq1nthrootDreg}
J_\mu \subset \liminf_{n\rightarrow\infty} J({\phi\circ q_n})\subseteq 
\limsup_{n\rightarrow\infty} K({\phi\circ q_n} )\subseteq \textup{Co}(S(\mu)),
\end{equation} 
if $J_\mu$ is Dirichlet regular. Moreover in both cases let $V_\epsi := \{z\in\C\; |\; g_{\Om_\mu}(z) \geq \epsi\}$.
Then 
\begin{equation}
\label{mainq2nthroot}
\forall\; \epsi>0\;:\qquad \lim_{n\to\infty} \Cpct(V_\epsi\cap K({\phi\circ q_n})) = 0 .
\end{equation}

\end{theorem}

\begin{corollary}
\label{cor:mainmonic}
Let $\phi$ be affine and let $\mu$ be $n$-th root regular with $\Cpct(S(\mu)) = 1$.
Then the sequence $\{P_n(z) = P_n(\mu; n)\}_n$ of monic orthogonal polynomials for $\mu$ 
satisfies 
\begin{equation}
\label{mainmonic1}
\overline{J_\mu \sm E} \subseteq 
\liminf_{n\rightarrow\infty} J({\phi\circ P_n} )\subseteq 
\limsup_{n\rightarrow\infty} K({\phi\circ P_n}) \subseteq \Co(S(\mu)).
\end{equation}   
Moreover let $V_\epsi := \{z\in\C\; |\; g_{\Om_\mu}(z) \geq \epsi\}$. 
Then 
\begin{equation}
\forall\; \epsi>0\;:\qquad \lim_{n\to\infty} \Cpct(V_\epsi\cap K({\phi\circ P_n})) = 0
\end{equation}
\end{corollary}

\begin{proof}
The monic sequence $P_n = p_n(\mu; z)/\gamma_n$ is asymptotically orthonormal for  $\mu$, 
so Theorem \ref{thm:mainq} implies the result. 
\end{proof}

\begin{corollary}
\label{cor:mainmonictwo}
Let $\mu$ be $n$-th root regular, let $\phi(z) = az + b$ with $|a| = 1/\Cpct(S(\mu)$ 
and let $\nu = \phi_\ast(\mu)$. 
Then the sequence $\{P_n(z) = P_n(\nu; n)\}_n$ of monic orthogonal polynomials for $\nu$ 
satisfies 
\begin{equation}
\label{mainmonic1}
\overline{J_\nu \sm E} \subseteq 
\liminf_{n\rightarrow\infty} J({P_n}) \subseteq 
\limsup_{n\rightarrow\infty} K({P_n}) \subseteq \Co(S(\nu)),
\end{equation}
where $E\subset S(\nu)$ is the set of Dirichlet non-regular points in $S(\nu)$.
Moreover let $V_\epsi := \{z\in\C\; |\; g_{\Om_\mu}(z) \geq \epsi\}$. 
Then 
\begin{equation}
\forall\; \epsi>0\;:\qquad \lim_{n\to\infty} \Cpct(V_\epsi\cap K({P_n})) = 0 .
\end{equation}
\end{corollary}

\begin{remark}
If there exists a subsequence $n_k$ such that $\lim_{k \rightarrow\infty} \vert a_{n_k} \vert^{\frac{1}{n_k}} = 1$, then Theorem \ref{thm:mainq} holds for the sequence $n_k$. 
\end{remark}

The final section of the paper is devoted to applications of the Main Theorem. It focuses on the limits of sequences of natural parameter spaces of polynomials coming from affine deformations of the sequences of orthonormal polynomials and  monic orthogonal polynomials for equilibrium measures. These applications are motivated by the  classical result \citep{barnsleyinvariantmeasures}. 

\begin{theorem}[\citep{barnsleyinvariantmeasures} Thm. 3]
\label{thm:barnsleyinvariantmeasures}
Let $f(z) $ be a monic centered polynomial of degree $d \geq 2$. Let $\mu$ denote the unique measure of maximal entropy for $f$, also known as the equilibrium measure for the Julia set $J(f)$. Then the monic orthogonal polynomials $\{ P_n (\mu ; z ) \}_n$ satisfy
\begin{align*}
P_1 (z) = z, \qquad P_{kd} (z) =P_k (f(z)), \qquad P_{d^k} (z) = f^k(z). 
\end{align*}
\end{theorem}
In particular, the sequence of iterates $\{ f^n \}_n$ is a subsequence of the sequence of monic orthogonal polynomials. The following corollary completes the discussion in \cite{lunde} on 
how to interpret the monic polynomials in the above theorem.
 
\begin{corollary}
\label{cor:seqofiterates}
Let $f$ and $\mu$ be as above. And let $\phi(z)$ be affine. Then
\begin{equation}
J(f) \subseteq \liminf_{n\rightarrow\infty} J({\phi\circ P_n}) \subseteq 
\limsup_{n\rightarrow\infty} K({\phi\circ P_n}) \subseteq \Co(J(f)).
\end{equation}   
Moreover let $V_\epsi := \{z\in\C\; |\; g_{\Om_\mu}(z) \geq \epsi\}$. 
Then 
\begin{equation}
\forall\; \epsi>0\;:\qquad \lim_{n\to\infty} \Cpct(V_\epsi\cap K({\phi\circ P_n})) = 0 .
\end{equation}
\end{corollary}
\begin{proof}
We have $S(\mu) = J(f)$, and $\Cpct (J(f)) = 1$ since $f$ is monic, see \eqref{eq:capacityleadingcoefficient}. 
Moreover, the equilibrium measure $\mu$ is $n$th root regular by the Erdös-Turán criterion \citep[Thm. 4.1.1]{stahltotik}.
 Thus 
Corollary \ref{cor:mainmonic} applies. 
\end{proof}

\section{Background}

Recall that $S(\mu)$  is the support of $\mu \in \mathfrak{B}$, 
where $\mathfrak{B}$ is the set of  Borel probability measures on $\C$ 
with compact nonpolar support. Moreover, recall that
$\Omega = \Omega_\mu$ denotes the unbounded connected component 
of $\C \sm S(\mu)$. 

\begin{theorem}[\citep{stahltotik} Lemma 1.1.3, originally proven by Fej{\'e}r \cite{fejer}]
\label{thm:stahl211}
All zeros of the orthonomal polynomial $p_n ( \mu ; z)$  are contained in the convex hull \ $\textup{Co}(S(\mu))$, and for any compact set $V \subseteq \Omega$ the number of zeros of $p_n (\mu ;  z)$ on $V$ is bounded as $n \longrightarrow \infty$.
\end{theorem}

In the following  $g_\Omega$ denotes  the Green's function for $\Omega$ with pole at infinity, while $g_\mu$ is the minimal-carrier Green's function. 
The minimal-carrier Green's function $g_\mu$ arises from an extension of the notion of Green's functions with pole at infinity to arbitrary Borel sets with bounded complement \citep[see][Appendix A.V]{stahltotik}.
 In particular, one can consider the unbounded connected component of  $\C \sm C$ 
 for $C$ a bounded carrier of $\mu$, i.e. $C$ a bounded measurable set satisfying $\mu (\C \sm C ) = 0$.  This unbounded component is not necessarily open. The Green's function $g_{\C \sm C}$ is then defined as  the unique non-negative subharmonic function, which is harmonic and positive on the interior of $\C \sm C $ and satisfies
\begin{equation}
g_{\C \sm C} (z) = \log \vv z \vv - \log \Cpct (C) + o(1),
\end{equation} 
which reflects the properties of the Green's function $g_\Omega$.
 The minimal-carrier Green's function $g_\mu$ is the supremum over all such $g_{\C \sm C}$ and  satisfies 
\begin{equation}\label{eq:minimalcarrier}
g_\mu (z) = \log \vv z \vv - \log c_\mu + o(1),
\end{equation} 
where $c_\mu$ is the minimal-carrier capacity, i.e. the infimum of the capacities of the bounded carriers of $\mu$.

\begin{theorem}[\citep{stahltotik} Thm. 1.1.4]
\label{thm:stahl114}
For a Borel probability measure $\mu$,
\begin{equation}
\label{eq:ST1}
\limsup_{n \rightarrow \infty} \left\vert p_n (\mu ; z) \right\vert^\frac{1}{n} \leq \eul^{g_\mu(z)}
\end{equation}
locally uniformly in $\mathbb{C}$, and 
\begin{equation}
\label{eq:ST2}
\liminf_{n \rightarrow \infty} \left\vert p_n (\mu ; z) \right\vert^{\frac{1}{n}} \geq \eul^{g_\Omega (z)}
\end{equation}
locally uniformly in $\mathbb{C} \setminus \textup{Co}(S(\mu))$.
In $\textup{Co}(S(\mu)) \cap \Omega$ the asymptotic lower bound holds true only in capacity, that is, for every compact set $V \subseteq \Omega$ and every $\epsi > 0$,
\begin{equation}
\label{eq:ST3}
\lim_{n \rightarrow \infty} \Cpct \left( \{ z \in V : \left\vert p_n (\mu ; z) \right\vert^{\frac{1}{n}} < \eul^{g_\Omega (z)} - \epsi \} \right) = 0 .
\end{equation}
In $\textup{Co}(S(\mu)) \cap \Omega$ the lower bound can also be given in the following form: \\ For every infinite subsequence $N \subseteq \mathbb{N}$,
\begin{equation}
\limsup_{n \rightarrow \infty, \; n \in N} \vert p_n (\mu ; z) \vert^{\frac{1}{n}} \geq \eul^{g_\Omega (z)}
\end{equation}
nearly everywhere in $\textup{Co}(S(\mu)) \cap \Omega$ , and on the outer boundary  $\partial \Omega$ of $S(\mu)$
\begin{equation}
\limsup_{n \rightarrow \infty, \; n \in N} \vert p_n (\mu ; z) \vert^{\frac{1}{n}}\geq 1
\end{equation}
nearly everywhere on $\partial \Omega$.
\end{theorem}

The relation \eqref{eq:ST1} holds locally uniformly in $\C$, which means that for any $z \in \C$ and $z_n \longrightarrow z$ we have 
$\limsup_{n \rightarrow \infty} \left\vert \phat_n (\mu ; z_n) \right\vert^\frac{1}{n} \leq \eul^{g_\mu(z)}$. Similarly for \eqref{eq:ST2}. 

We have the following corollary.

\begin{corollary}
\label{cor:logST}
For a Borel probability measure $\mu$,
\begin{equation}
\label{eq:logST1}
\limsup_{n \rightarrow \infty} \frac{1}{n} \log^+ \left\vert p_n (\mu ; z) \right\vert \leq {g_\mu(z)}
\end{equation}
locally uniformly in $\mathbb{C}$, and 
\begin{equation}
\label{eq:logST2}
\liminf_{n \rightarrow \infty} \frac{1}{n} \log^+ \left\vert p_n (\mu ; z) \right\vert \geq {g_\Omega (z)}
\end{equation}
locally uniformly in $\mathbb{C} \setminus \textup{Co}(S(\mu))$.
In $\textup{Co}(S(\mu)) \cap \Omega$ the asymptotic lower bound holds true only in capacity, that is, for every compact set $V \subseteq \Omega$ and every $\epsi > 0$,
\begin{equation}
\label{eq:logST3}
\lim_{n \rightarrow \infty} \Cpct \left( \{ z \in V : \frac{1}{n}\log^+ \left\vert p_n (\mu ; z) \right\vert <{g_\Omega (z)} - \epsi \} \right) = 0 .
\end{equation}
\end{corollary}

\begin{proof}
The relations \eqref{eq:logST1} and \eqref{eq:logST2} follows directly from the corresponding relations \eqref{eq:ST1} and \eqref{eq:ST2}.
For the last statement consider $ \{ z \in V : \frac{1}{n}\log^+ \left\vert p_n (\mu ; z) \right\vert <{g_\Omega (z)} - \epsi \} = \{ z \in V :  \left\vert p_n (\mu ; z) \right\vert^ \frac{1}{n} <\eul^{{g_\Omega (z)} - \epsi} \}$. 
Since $\eul^{{g_\Omega (z)} - \epsi} = \eul^{g_\Omega (z)} -  \eul^{g_\Omega (z)} (1 - \eul^{-\epsi} ) <  \eul^{g_\Omega (z)} -  (1 - \eul^{-\epsi} )$ as $g_\Omega (z) \geq 0$, it follows that $ \{ z \in V : \frac{1}{n}\log^+ \left\vert p_n (\mu ; z) \right\vert <{g_\Omega (z)} - \epsi \} \subseteq \{ z \in V :  \left\vert p_n (\mu ; z) \right\vert^ \frac{1}{n} <\eul^{g_\Omega (z)} - \epsi' \}$ with $\epsi ' = 1 - \eul^{-\epsi}$.  Now \eqref{eq:ST3} implies \eqref{eq:logST3}. 
\end{proof}

It is customary in the theory of orthogonal polynomials to require that the leading coefficient $\gamma_n$  in \eqref{eq:sequenceelement} satisfies $\gamma_n > 0$, which gives uniqueness of the sequence in $L^2(\mu)$. 
Uniqueness aside, this requirement is often superfluous when proving statements on orthonormal polynomials.
One can define  $\gammahat_n := \epsi_n \gamma_n \in \C \sm \{ 0 \}$ for any 
sequence $\{\epsi_n\}_n$ with $\vv \epsi_n \vv = 1 $ 
and so obtain $\{  \phat_n (\mu ; z ) \}_n := \{ \epsi_n p_n (\mu ; z ) \}_n$, 
another orthonormal sequence for $\mu$. 
Since $\vv\vv p_n \vv\vv_{L^2 (\mu)} =  \vv\vv \epsi_n p_n \vv\vv_{L^2 (\mu)} $, 
such orthonormal sequences are not uniquely determined by the $L^2 (\mu)$ norm. 
Varying the argument of $\gammahat_n$ yields different sequences of orthonormal polynomials, while fixing the sequence $\{\epsi_n\}_n$ gives the unique orthonormal sequence with leading coefficient $\gammahat_n$. 
We shall refer to any sequence $\{  \phat_n (\mu ; z ) \}_n$ as an orthonormal sequence 
and to the sequence $\{ p_n (\mu ; z ) \}_n$ as the canonical orthonormal sequence, 
when need be.

Since $|\phat_n (\mu ; z )| = |p_n(\mu;z)|$ for each $n$ and the two sets of zeros for the polynomials are identical, the classical result \citep[Thm. 1.1.4]{stahltotik} 
is invariant under this modification:

\begin{proposition}
Replacing  $p_n (\mu ; z)$ with $\phat_n (\mu ; z )$ preserves all inequalities in Theorem \ref{thm:stahl114}.
\end{proposition}

More is possible. We also consider orthogonal polynomial sequences which 
asymptotically are orthonormal as described in Definition \ref{def:asymptoticallyorthonormal}.

If $\vv a_n \vv = 1$ for all $n$, then $\{q_n (\mu ; z)\}_n$ 
is an orthonormal sequence for $\mu$ 
as in the above discussion and the above proposition applies. 
In general the slightly weaker Corollary \ref{cor:logST} of \citep[Thm. 1.1.4]{stahltotik} 
holds for asymptotically orthonormal sequences $\{q_n(\mu ; z)\}_n$, 
this version however suffices for the results of this present paper.

\begin{proposition}
\label{prop:invarianceSTq}
Replacing $p_n (\mu ; z)$ with an asymptotically orthonormal sequence $q_n (\mu ; z)$ preserves all inequalities in Corollary \ref{cor:logST}.
\end{proposition}

\begin{proof}
Since $a_n \neq 0$ the zeros of $\qhat_n(\mu ; z)$ and $p_n(\mu ; z)$ coincide.
Moreover, $\lim_{n \rightarrow \infty}\frac{1}{n}\log^+ \vv a_n \vv = 0$, thus the relations \eqref{eq:logST1} and \eqref{eq:logST2} hold for $\qhat_n (\mu ; z ) = a_n \phat_n (\mu ; z)$ as well.  Regarding \eqref{eq:logST3}, choose $N$ such that $\left\vert \frac{1}{n} \log \vert a_n \vert \right\vert < \frac{\epsi}{2}$ for $n > N$. Then $ \{ z \in V : \frac{1}{n}\log^+ \vv \qhat_n (\mu ; z ) \vv < g_\Omega (z) - \epsi \} \subseteq  \{ z \in V : \frac{1}{n}\log^+ \vv p_n (\mu ; z ) \vv < g_\Omega (z) - \frac{\epsi}{2} \} $ for $n > N$, so \eqref{eq:logST3} implies the result for $\qhat_n(\mu ; z)$.
\end{proof}

\begin{remark}
\label{rmk:invarianceSTq}
Note that also \eqref{eq:ST1} and \eqref{eq:ST2} hold for any asymptotically orthonormal sequence.
\end{remark}

Our results are phrased in terms of the Hausdorff distance between compact subsets of $\C$. 
The space of non-empty compact subsets of $\C$ equipped with the Hausdorff distance is a complete metric space \cite[Chapter II]{barnsley}.
 For  a  bounded sequence $\{ S_n \}_n$ of non-empty compact sets define the sets
\begin{align} 
\liminf_{n\rightarrow\infty} S_n &= \{ z \in \C \mid \exists \{ z_n \}, S_n \ni z_n \underset{n\rightarrow\infty}  {\longrightarrow } z \},
\\
\limsup_{n\rightarrow\infty} S_n  &= \{ z \in \C \mid \exists \{ n_k \}, n_k \nearrow \infty \text{ and } \exists  \{ z_{n_k} \}, S_{n_k} \ni z_{n_k} \underset{n\rightarrow\infty}  {\longrightarrow } z \},
\end{align}
 which are both compact \citep[Lemma 3.1]{lunde}. The set $\displaystyle\liminf_{n\rightarrow\infty} S_n $ might be empty however.

The notions of $\liminf$ and $\limsup$ apply to the sequences of Julia sets and filled Julia sets. Both the Julia set and filled Julia set are compact and moreover the sequences 
$\{ J({q_n}) \}_n$ and $\{ K({q_n}) \}_n$ are bounded for $n$ sufficiently large. This follows from \citep[Lemma 2.2]{lunde}, which by Proposition \ref{prop:invarianceSTq} and inspection holds for $q_n$ as well as $p_n$, see also Proposition \ref{prop:lundeinvarianceq} below.

We shall moreover be interested in closed bounded subsets of $\Cstar\times\C$. 
The above notion of Hausdorff distance between compact sets trivially generalises to any metric space. However we may also extend the notion of Hausdorff distance to bounded 
relatively closed subsets of $\Cstar\times\C$ by defining for any two non-empty 
bounded and relatively closed subsets $A, B \subset\Cstar\times\C$ 
$$
d_H(A, B) = d_H(\overline{A}, \overline{B})
$$
where closures are in the ambient $\C^2$ and the right hand $d_H (\cdot , \cdot)$ denotes the Euclidean Hausdorff distance between compact sets.

Finally, the subset $\mathfrak{B}_+ \subseteq \mathfrak{B}$ 
here defined as
\begin{equation}
\label{eq:Bplus}
\mathfrak{B}_+ = \{ \mu \in \mathfrak{B} \mid \liminf_{n\rightarrow\infty} \Cpct (K({p_n (\mu ; z)})) > 0 \}
\end{equation}
 is defined in terms of the leading coefficient $\gamma_n$ in \citep{lunde}. 
In particular $ \mathfrak{B}_+ := \{ \mu \in \mathfrak{B} \mid \limsup_{n\rightarrow\infty} \gamma_n^{\frac{1}{n}} < \infty \}$. Note that it follows directly from this definition and \eqref{eq:nthrootreg1}  that  if $\mu \in \mathfrak{B}$ is $n$th root regular, then $\mu \in \mathfrak{B}_+$. 
 The equivalence of the two descriptions easily 
follows from the relation 
\begin{equation}
\label{eq:capacityleadingcoefficient}
\Cpct (K({p_n (\mu ; z)}) )= \frac{1}{\gamma_n^{\frac{1}{n-1}}} ,
\end{equation}
see \citep[Thm. 6.5.1]{ransford}.
  We choose the above description under the impression that it is the more illustrative of the geometric properties relevant here.

\section{Conformal Invariance}
\label{sec:conformalinvariance}

Given  \eqref{eq:polynomialtransformation} we can verify that the above relations from \citep{stahltotik} are invariant under conformal equivalence of measures.
First recall that $\nu = \nu_{a,b} := \phi^{-1}_\ast (\mu) = \phi^\ast (\mu)$, 
where $\phi(z) = \aphi(z) = az+ b$, $a\in\C\setminus\{0\}$, $b\in \C$.
Moreover 
\begin{equation}
\label{eq:supporttransformation}
S(\mu) = \phi (S (\nu))
\end{equation}
and $\{\phat_n (\nu ;  z):= \phat_n (\mu ; \phi(z))\}_n$ is an 
orthonormal sequence for  $\nu$, whenever 
$\{\phat_n (\mu ;  z)\}_n$ is an orthonormal sequence for $\mu$. 
And similarly $\{q_n (\nu ;  z):= q_n (\mu ; \phi(z))\}_n$ 
is an asymptotically orthonormal sequence for $\nu$, whenever 
$\{q_n (\mu ;  z) = a_np_n(\mu; z)\}_n$ is an asymptotically orthonormal sequence for $\mu$.

\begin{proposition}
\label{prop:conformalinvariance}
Replacing $\mu$ with $\nu$ and $\phat_n (\mu ; z) $  with $\phat_n (\nu ;  z)$ 
preserves all inequalities in Theorem \ref{thm:stahl114} 
as well as the relative position of the zeros and the support. 
\end{proposition}

\begin{proof}
Assume $\phat_n (\nu ; z_0 ) = 0$, so that $z_0 \in \Co(S (\nu))$. Then $\phi (z_0) \in \Co( S(\mu)) $ by \eqref{eq:supporttransformation} and  $\phat_n (\mu ; \phi (z_0) ) = 0$ by \eqref{eq:polynomialtransformation}. Hence the position of the zeros relative to the support is preserved. 

Moreover, we have that
\begin{equation}
\frac{1}{n} \log \vv \phat_n (\nu ; z  ) \vv = \frac{1}{n}  \log \vv \phat_n (\mu ; \phi (z)  ) \vv
\end{equation}
and as Green's functions are conformally invariant  \citep[see][Thm. 4.4.4]{ransford} these  transform in the same way, i.e. 
\begin{equation}
\label{eq:greensfctconformalinvariance}
g_{\Omega ( \nu )}(z) =  g_{\Omega (\mu )}(\phi(z)).
\end{equation}
So the relations between $\frac{1}{n} \log \vert \phat_n ( \nu ; z ) \vert$ and $g_{\Omega ( \nu )} (z)$ are the same as the relations between $\frac{1}{n} \log \vert \phat_n ( \mu  ;\phi ( z)) \vert$ and $g_{\Omega(\mu)} (\phi (z))$. Similarly for the relations involving the minimal carrier Green's functions $g_\nu (z)$ and $g_\mu (\phi (z ))$. Hence the asymptotic bounds are preserved under a conformal pullback of the measure as well. 
\end{proof}

We have a similar result for $q_n (\mu ; z)$ on the invariance of the weaker version of the results in \cite{stahltotik} given here as Corollary \ref{cor:logST}. 
The proof is as of the above proposition.

\begin{proposition}
\label{prop:conformalinvarianceq}
Replacing $\mu$ with $\nu$ and $q_n (\mu ; z) $  with $q_n (\nu ;  z)$ preserves all inequalities in Corollary \ref{cor:logST} as well as the relative position of the zeros and the support. 
\end{proposition}

The above propositions indicate that from the point of view of orthogonal polynomials 
there is no difference between working with $\mu$ or any equivalent measure 
$( \aphi) ^\ast (\mu)$. 
We therefore say that asymptotically orthonormal polynomials are conformally invariant. However, as was noted in the introduction the dynamics of orthogonal polynomials 
$q_n(\mu, z)$ and $q_n(\nu, z)$ in general will be quite different. 
See for example Section \ref{subsec:multibrot}.

Recall that two polynomials $p$ and $q$ are conformally conjugate if there is a biholomorphic hence affine map $\phi$  such that $q \circ \phi = \phi \circ p$, or equivalently
\begin{equation}
\label{eq:confequivdynsys}
p = \phi \circ q \circ \phi^{-1} .
\end{equation}
If $p$ and $q$ are conformally conjugate under $\phi$, then
\begin{equation}
\label{eq:equivdynamics}  
K(p) = \phi(K(q) ), \quad J(p) = \phi(J(q) ) .
\end{equation}

 Based on the above propositions one could think  that $q_n (\mu ; z)$ and $q_n (\nu ; z)$ 
are conformally equivalent dynamical systems, however it  is evident from \eqref{eq:polynomialtransformationq} that this is not true in general. 
We will use this observation to describe affine deformations of $q_n (\mu ; z)$. 

The  dynamics of the sequence $ \{ p_n (\mu ; z) \}_n$ is described in the paper  ``Julia Sets of Orthogonal Polynomials'' \citep{lunde}.
 In \citep{lunde}  a relationship is presented between dynamical properties of the sequence $\{ p_n  (  \mu ; z ) \}_n$ of orthonormal polynomials and the compact support $S(\mu)$ of $\mu \in \mathfrak{B}_+$. In particular that
\begin{equation}
\label{eq:lunde}
\overline{J _ \mu \sm E_{\mu}} \subseteq \liminf_{n\rightarrow\infty} J ( {p_n (\mu ; z )}  )\subseteq \limsup_{n\rightarrow\infty} K( {p_n (\mu ; z )}) \subseteq \textup{Co}(S(\mu)) ,
\end{equation}
where $E_\mu$ is the exceptional set for the minimal carrier Green's function $g_\mu$, and 
\begin{equation}
\label{eq:lundecapacitet}
\forall\; \epsi>0\;:\qquad \lim_{n\to\infty} \Cpct(V_\epsi\cap K({ p_n(\mu ; z )})) = 0 ,
\end{equation}
where $V_\epsi := \{z\in\C\; |\; g_{\Om_\mu} (z) \geq \epsi\}$.
When $\mu$ is $n$th root regular $E_\mu$ is the exceptional set for the Green's function $g_\Omega$ \citep[see][Thm. 1.3 and Corollary 4.2]{lunde}.

The results in \citep{lunde} are proven for the orthonormal $p_n (\mu ; z)$ with $\gamma_n > 0$, but we can replace the sequence $\{p_n (\mu ; z)\}_n$ with any asymptotically orthonormal sequence $\{q_n (\mu ; z)\}_n$.
  
  \begin{proposition}
  \label{prop:lundeinvarianceq}
Let $\mu\in\mathfrak{B}_+$ and let $\{q_n (z) = a_np_n(\mu; z)\}_n$ be 
any  asymptotically orthonormal sequence for $\mu$. Then 
\begin{equation}
\label{eq:lundeq}
\overline{J_\mu \sm E_{\mu}} \subseteq \liminf_{n\rightarrow\infty} J({q_n ( z )}) \subseteq \limsup_{n\rightarrow\infty} K({q_n ( z )}) \subseteq \textup{Co}(S(\mu)) ,
\end{equation}
where $E_\mu$ is the exceptional set for the minimal carrier Green's function $g_\mu$. Moreover, $E_\mu$ can be replaced by the exceptional set $E$ for the Green's function $g_{\Omega_\mu}$, when $\mu$ is $n$-th root regular.

Furthermore, let $V_\epsi := \{z\in\C\; |\; g_{\Om_\mu}(z) \geq \epsi\}$. 
Then 
\begin{equation}
\forall\; \epsi>0\;:\qquad \lim_{n\to\infty} \Cpct(V_\epsi\cap K({ q_n (z)})) = 0 .
\end{equation}
\end{proposition}
  
\begin{proof}
Inspection shows that Proposition \ref{prop:invarianceSTq} implies that the results in \citep{lunde} remain valid for $q_n (\mu ; z)$. In particular we observe the following. 
  
Let $K_\mu = \C \sm \Omega_\mu$ denote the filled support of $\mu \in \mathfrak{B}$ 
and $g_n$ 
the Green's function for $\C \sm K( q_n ( z))$.
If $K_\mu \subseteq \mathbb{D}(0,R)$, then there is $N$ such that 
$K( q_n ( z)) \subset q_n^{ \ -1} ( \overline{\mathbb{D} (0,R)}) \subset  \mathbb{D} (0 ,R)$ for $n > N$. This containment is proven for $p_n (\mu ; z)$ in \citep[Lemma 2.2]{lunde}. But the only elements particular to orthonormality in the proof  is \eqref{eq:logST2} and the position of the zeros, so by Proposition \ref{prop:invarianceSTq} the same is true for $q_n ( z)$. 
With this result we directly obtain the main technical tool 
\citep[Prop. 2.3]{lunde}, 
that for any $\mu \in \mathfrak{B}_+$ there is a natural number $N$ and a constant $M > 0$ such that  
  \begin{equation}
  \label{eq:CHPPmaintoolq}
  \left\vert\left\vert g_n (z) - \frac{1}{n} \log^+\vert \qhat_n (\mu ; z) \vert \right\vert\right\vert_\infty \leq \frac{M}{n}.
  \end{equation}
As in \cite[Prop. 2.7]{lunde} this allows us to relate the dynamical Green's functions $g_n$ to the Green's functions $g_\mu$ and $g_\Omega$ associated to $\mu$ through Corollary \ref{cor:logST} which by Proposition \ref{prop:invarianceSTq} 
holds for $q_n(\mu ; z)$ as well. 
An inspection of the proofs of \eqref{eq:lunde} and \eqref{eq:lundecapacitet} \citep[Thm. 1.3] {lunde} 
reveals that with the above results preserved we can replace 
$p_n (\mu ; z)$ with $q_n (\mu ; z)$ in \eqref{eq:lunde} and \eqref{eq:lundecapacitet}.
\end{proof}

In the rest of this paper we explore the measure dependency of this geometric relationship, i.e. how a conformal transformation $\mu\mapsto \nu_{a,b} = (\aphi)^\ast(\mu)$ 
of the measure affects the relative position of the dynamical sets and the support.
The dynamical systems given by the polynomials are not conformally conjugate, which indicates that the relative position is not preserved under conformal equivalence of measures. 
While the conclusions in the above theorems from \citep{stahltotik} are invariant under conformal equivalence of measures, the conclusions in \citep{lunde} have a new content for every nontrivial pullback of a measure $\mu \in \mathfrak{B}_+$. 
Theorem \ref{thm:mainq} combined with the results in \citep{lunde} and Proposition \ref{prop:lundeinvarianceq} shows that even though the dynamics is dependent on the normalization of the measure some structure is still retained.

\section{Proof of Theorems \ref{thm:mainq} and \ref{thm:mainregular}}

As remarked above, the polynomials in an orthogonal sequence $\{q_n (\mu ; z)\}_n$ 
of polynomials for $\mu$ are not conformally conjugate to the polynomials in the corresponding orthogonal sequence $\{q_n (\nu ; z) = q_n(\mu ; \phi(z)\}_n$ for the pullback measure $\nu$, when viewed as dynamical systems. 
The key for proving Theorem \ref{thm:mainq} is the following simple observation.

\begin{lemma}
\label{prop:conformalconjugate}
Let $\mu \in \mathfrak{B}$, and let $\phi$ be affine.
Then for every $n\geq 1$ the polynomials $q_n(\nu ; z)$ and $\phi \circ q_n ( \mu ; z)$ 
are conjugate under $\phi$.
\end{lemma}
One can visualize this claim by the following commuting diagram,
\begin{center}
\begin{tikzpicture}[baseline= (a).base]
\node[scale=1.3] (a) at (0,0){
\begin{tikzcd}
\C \arrow{r}{   \phi  \circ q_n (\mu ; z) }  \arrow[swap]{d}{\phi} & [3em] \C  \arrow{d}{\phi} \\ 
\C \arrow[swap]{r}{     q_n (\nu ; z)   } & \C
\end{tikzcd}
};
\end{tikzpicture}
\end{center}

\begin{proof}
Since  $\phi\circ q_n(\mu ; z) = \phi\circ q_n(\mu ; \phi\circ\phi^{-1}(z)) = 
\phi\circ q_n(\nu ; \cdot)\circ\phi^{-1}(z)$ by \eqref{eq:polynomialtransformationq}  it follows from \eqref{eq:confequivdynsys} that $\phi \circ q_n (\mu ; z )$ and $q_n (\nu ; z)$    are conformally conjugate under $\phi$. 
\end{proof}

The results in \cite{lunde}, in particular \eqref{eq:lunde} and \eqref{eq:lundecapacitet} above, directly applies to the pullback measure $\nu = \phi^\ast (\mu)$ and its support $S(\nu)$. By Proposition \ref{prop:lundeinvarianceq} this relationship also holds for $q_n (\mu ; z)$. 

We  can now prove our main result, which for every affine map $\phi$
relates the support $S(\mu)$ to the sequences  $\{ J ( {\phi \circ q_n ( \mu ; z)} ) \}_n$ 
and $\{ K({\phi \circ q_n ( \mu ; z)})\}_n$ of Julia sets and filled Julia sets.

\begin{proof}[Proof of Theorems \ref{thm:mainq} and \ref{thm:mainregular}] 
Let $\mu \in \mathfrak{B}$, let $\phi$ be affine and let $\nu =  \phi ^\ast (\mu )$. 
Suppose $q _n(z) = q_n( \mu ; z) = a_n p_n (\mu ; z)$ is an orthogonal polynomial sequence for $\mu$. Then 
$$
a_np_n( \nu ; z) = q_n( \nu ; z) = q_n \circ \phi(z)
$$ 
is the corresponding  orthogonal polynomial sequence for $\nu$. 
Thus $S(\mu)$ and the Julia sets $J({\phi\circ q_n})$ 
and filled Julia sets $K({\phi\circ q_n})$ have the same relative position as 
the support $S(\mu) = \phi ( S(\nu) ) $ of the measure $\nu$ 
and the Julia sets $J({q_n\circ\phi})$ and filled Julia sets $K( {q_n\circ\phi})$,
since it follows from Lemma \ref{prop:conformalconjugate} 
and \eqref{eq:equivdynamics} that 
\begin{equation}
\begin{aligned}
\label{eq:affinemapsofjuliasets}
J({\phi \circ q_n ( z  )}) & = \phi (J({q_n \circ \phi(z) }) ) = \phi (J({q_n ( \nu ;z )}) ) 
\\
 K({\phi \circ q_n ( z  )})  & =   \phi (K({q_n \circ \phi(z) } )) =    \phi (K({q_n (\nu ; z )} )).
\end{aligned}
\end{equation}

Moreover,  if $\mu \in \mathfrak{B}_+$ and 
the sequence $\{q_n(z)\}_n$ is  asymptotically orthonormal for $\mu$, then
according to 
Proposition \ref{prop:lundeinvarianceq}
\begin{equation}
\label{main2}
\overline{J_\nu \sm E_{\nu}} \subseteq \liminf_{n\rightarrow\infty} J ({q_n (\nu ; z )} )
\subseteq \limsup_{n\rightarrow\infty} K({q_n (\nu ; z )}) \subseteq \textup{Co}(S(\nu)) ,
\end{equation}
where $E_\nu$ is the exceptional set for the minimal carrier Green's function 
$g_\nu$ for $\nu$. If moreover $\mu$ and hence $\nu$ is $n$th root regular 
then we can replace $E_\nu$ by the set $E$ consisting of Dirichlet non-regular 
outer boundary points of $S(\nu)$.

By \eqref{eq:affinemapsofjuliasets} this is equivalent to \eqref{mainq1} 
and to \eqref{mainq1nthroot} if $\nu$ is $n$th root regular as the relative positions of the sets are preserved.

Moreover, according to 
Proposition \ref{prop:lundeinvarianceq}
\begin{equation}
\forall\; \epsi>0\;:\qquad \lim_{n\to\infty} \Cpct(V_\epsi (\nu )\cap K({ q_n (\nu ; z)}) )= 0 
\end{equation}
with  $V_\epsi (\nu) = \{z\in\C\; |\; g_{ \Om _\nu }(z) \geq \epsi\} $. 
Since
\begin{align*}
\phi ( V_\epsi (\nu) \cap K({ q_n (\nu ; z)}) ) & = \phi (V_\epsi (\nu) ) \cap K({ \phi \circ q_n  (z)} ) 
\\ & = \{ z \mid g_{\Om _ \nu} (\phi^{-1} (z)) \geq \epsi \}\cap K({ \phi \circ q_n (z)} )
\\ & =  
V_\eps (\mu) \cap K (\phi \circ q_n (z))
\end{align*}
by \eqref{eq:greensfctconformalinvariance}, it follows that 
$
\Cpct  (  V_\eps (\mu ) \cap K( { \phi \circ q_n (z)} ) ) = a   \Cpct(V_\epsi (\nu )\cap K({ q_n (\nu ; z)}) )
$
by translation invariance of capacity. This proves \eqref{mainq2}  and in the case $\nu$ is $n$th root regular  \eqref{mainq2nthroot} and thus  
completes the proof of Theorems \ref{thm:mainq} and \ref{thm:mainregular}.
\end{proof}

\section{Applications}  
\label{sec:examples}

Let $K$ be a compact continuum with connected complement $\Omega = \C\setminus K$. 
Then $K$ is non-polar and its boundary $J$ is Dirichlet regular, being connected. 
Let $\mu$ denote the equilibrium measure on $K$. 
Then the support  $S(\mu)$ equals $J$.
Such a measure $\mu$ is an example of a measure for which the associated orthonormal polynomials $p_n ( \mu ; z)$ have regular $n$th root asymptotic behavior by the Erdös-Turán criterion \citep[Thm. 4.1.1]{stahltotik}.  
  
If we assume also that $K$ is convex we have the following special case of  Theorem \ref{thm:mainregular} and in particular of 
\eqref{mainq1nthrootDreg}:
\begin{equation}
\label{eq:mainspecial}
J   \subseteq \liminf_{n\rightarrow\infty} J({\phi \circ q_n (  z  )} )\subseteq \limsup_{n\rightarrow\infty} K({\phi \circ q_n (  z  )}) \subseteq  K 
\end{equation} 
for every affine map $\phi$ and for any  asymptotically orthonormal sequence 
$\{q_n(z) = a_np_n( \mu ; z )\}$ for $\mu$.

Moreover if $\Cpct(K) = 1$ then according to Corollary \ref{cor:mainmonic}
\begin{equation}
\label{eq:mainspecialcpctone}
J    \subseteq \liminf_{n\rightarrow\infty} J({\phi \circ P_n ( \mu ;  z  )} ) \subseteq \limsup_{n\rightarrow\infty} K({\phi \circ P_n ( \mu ; z  )}) \subseteq  K .
\end{equation} 
where $P_n( \mu ; z )$ is the $n$th monic orthogonal polynomial for $\mu$.

In addition to the position of sequences of Julia sets and filled Julia sets, we study the connectedness of these sets.  Recall that the Julia set  for a polynomial is connected if and only if all finite critical points are contained in the filled Julia set. On the other hand, if all critical points escape to infinity under iteration, the Julia set is totally disconnected; see, e.g., \cite{milnor}.
We consider $\aphi \circ q_n$ as a two-parameter family and describe its connectedness locus, i.e. the set of parameter values $(a,b)$ such that $J({\aphi \circ q_n})$ is connected.

\subsection{The family $\mathbf{z^n + c}$}
\label{subsec:multibrot}

Let $\mu \in \mathfrak{B}_+$ denote the standard Lebesgue probability measure on the unit circle $\mathbb{S}^1$, i.e. the arc length measure on the unit circle normalized so the circle has arc length 1. Thus $S(\mu) = \mathbb{S}^1$ and $\rm{Co}(S (\mu)) = \overline{\D} = K (\mu)$. Then  $\mu$ is the equilibrium measure on $\mathbb{S}^1$, so $\mu$ has regular $n$th root asymptotic behavior.
It is well-known that $p_n (\mu ; z) = z^n$ is the unique orthonormal sequence for $\mu$.

Now let $\phi_c(z)=z + c$ for $c \in \C$ and let $\nu =( \phi_c )^\ast (\mu)$. By conformal equivalence, $\nu$ is $n$th root regular with Dirichlet regular support as well.
 It follows from  Lemma \ref{prop:conformalconjugate} that
$ p_n (\nu ; z  )$ is conformally conjugate to $\phi_c\circ p (\mu ; z ) = z^n + c$ and as the filled support $K(\mu) = \overline{\mathbb{D}}$ is  convex  \eqref{eq:mainspecial} implies that
\begin{align}
\label{ex2}
\mathbb{S}^1  
\subseteq \liminf_{n\rightarrow\infty} J({z^n +c} ) \subseteq \limsup_{n\rightarrow\infty} K({z^n +c}) \subseteq  
  \overline{\D}.
\end{align}

We compare this to results in \citep{boyd} and \citep{kaschner15}.

\begin{theorem}[\citep{boyd} Thm. 1.2]
Let $c \in \C$.\\
If $c \in \C \sm \overline{\D}$, then
\begin{align*}
\lim_{n\rightarrow\infty} J({z^n + c}) = \lim_{n\rightarrow\infty} K({z^n + c}) = \mathbb{S}^1.
\end{align*}
If $c \in \D$, then 
\begin{align*}
\lim_{n\rightarrow\infty} J({z^n + c}) = \mathbb{S}^1 \text { and } \lim_{n\rightarrow\infty} K({z^n + c}) =  \overline{\mathbb{D}}.
\end{align*}
If $c \in \mathbb{S}^1$, then if $\displaystyle\lim_{n\rightarrow\infty} J({z^n + c})$ and/or $ \displaystyle\lim_{n\rightarrow\infty} K({z^n + c})$ (and/or any $\liminf$ or $\limsup$) exists, it is contained in $\overline{\D}$. %
\end{theorem}

Notice that \eqref{ex2} holds for all $c \in \C$, so in particular for $c \in \mathbb{S}^1$.
We therefore  have the following corollary of Theorem \ref{thm:mainq},  which extends \citep[Lemma 3.7]{boyd} \citep[Lemma 3]{kaschner15}.

\begin{corollary} 
If $c \in \mathbb{S}^1$, then 
\begin{align*}
\mathbb{S}^1 \subseteq \liminf_{n \rightarrow \infty} J({z^n + c} ) \subseteq \limsup_{n\rightarrow\infty} K({z^n +c}) \subseteq   \overline{\D}.
\end{align*}
In particular,   if $ \{ J({z^{n_k} + c}) \}_k$ or $\{ K({z^{n_k} + c} )\}_k$ is a convergent subsequence, then  $\displaystyle\lim_{k\rightarrow\infty} J({z^{n_k} + c})$ respectively $ \displaystyle\lim_{k\rightarrow\infty} K({z^{n_k} + c})$ contains the unit circle $\mathbb{S}^1$. 
\end{corollary}

\begin{remark}
The connectedness locus $\mathcal{C} ( z^n + c) =  \mathcal{M}_n $ approaches the closed disc as $n$ tends to infinity  \citep[Thm. 1.1]{boyd}. Notice that by the above the only critical value $c$ belongs to $\C \sm \limsup_{n\rightarrow\infty} K({z^n +c}  )$ if $\vert c \vert > 1$. We thus obtain
$ \lim_{n\rightarrow\infty} \mathcal{M}_n \subseteq \overline{\mathbb{D}}$.
\end{remark}

\subsection{The Chebyshev Polynomials}

The Julia set for $p_2 (z) = z^2 - 2 $ is the real interval $[-2,2]$ \cite[see][Section 7]{milnor}.
The polynomial $p_2 (z)$ is the second Chebyshev polynomial in the sequence defined by
\begin{equation}
p_{n+1} (z) = z p_n (z) - p_{n-1} (z)
\end{equation}
with $p_1 (z) = z$ and $p_2 (z) $ as above. In fact, any Chebyshev polynomial $p_n (z)$ with $n \geq 2$ has Julia set $[-2,2]$.

The sequence of Chebyshev polynomials $p_n (z) = p_n (\mu ; z)$ is the sequence of monic orthogonal polynomials for the equilibrium measure $\mu$ on $[-2,2]$. This real interval satisfies that $\Cpct([-2,2]) = 1$ \cite[Cor. 5.2.4]{ransford}, thus Corollary \ref{cor:mainmonic} applies. 

Combining with \eqref{eq:mainspecialcpctone} it follows 
that 
\begin{equation}
[-2, 2] \subseteq \liminf_{n\rightarrow\infty} J({\aphi \circ p_n ( z)}) \subseteq \limsup_{n\rightarrow\infty} K({\aphi \circ p_n ( z)}) \subseteq [-2,2].
\end{equation}
In particular, for any $a \in \C \sm \{ 0 \}$ and $b \in \C$ the sequence $\aphi\circ p_n$ has Julia sets which converge to $[-2,2]$,
\begin{equation}
\label{eq:chebyshevlimits}
\lim_{n \rightarrow \infty} J({\aphi\circ p_n ( z)})  =   \lim_{n \rightarrow \infty} K({\aphi\circ p_n ( z)})
= [-2,2] .
\end{equation}

\subsubsection{Connectedness Loci}

If  $n \geq 3$, then the Chebyshev polynomial  $p_n (z)$ has $n-1$ critical points, but only two critical values, $z = \pm 2$. 

Consider the translation $\phi_{b} (z) = z + b$. The critical values of $\phi_{b} \circ p_n (z) = p_n (z) + b$ are therefore $b \pm 2$. Thus for any nontrivial translation at least one critical value is in $\mathbb{C} \sm [-2,2]$. So by \eqref{eq:chebyshevlimits} there is $N$ so that for $n\geq N$ at least one critical value is in $\C \sm K({p_n + b})$, which implies that $K({p_n + b})$ is disconnected. Thus the limit of the connectedness loci $\mathcal{C} ( p_n + b)$ of translations of the Chebyshev polynomials is 
\begin{align*}
\lim_{n\rightarrow \infty}\mathcal{C} ( p_n + b) =  \{ 0 \} .
\end{align*} 
Combining with \eqref{eq:chebyshevlimits} it follows that for any nontrivial translation of the sequence of Chebyshev polynomials $p_n$ and any $\epsi > 0$ there is $N$ so that for $n \geq N$ the Hausdorff distance $d_H (K({ p_n + b} ) , [-2,2] ) < \epsi$, while $K({ p_n + b })$ has uncontably many connected components.

Now consider the dilations $\phi_a \circ p_n (z) = a p_n (z)$, which have critical values $\pm 2 a$.  Recall that $a \in \C \sm\{ 0 \}$. 
If $a \notin [-1,1]$ there is $N$ so that for $n \geq N$ the filled Julia set $K({ap_n (z)})$ contains no critical values and is therefore a Cantor set.

This observation combined with the symmetry of the critical values allows us to determine the connectedness locus of the two-parameter family $ap_n + b$.  

\begin{theorem}
Let  $p_n (z)$ be the $n$-th  Chebyshev polynomial and let $\aphi (z) = az +b$. Then 
\begin{align}
\label{eq:chebyshevconnectedness}
\lim_{n\rightarrow \infty}\mathcal{C} ( a p_n + b) = \{ 
(a , b ) \in \mathbb{R}^2  \mid \ 0 <  \vert a \vert \leq 1 ,  \vert b \vert  \leq 2 ( 1 - \vert a \vert) \
 \}.
\end{align}
\end{theorem}

\begin{proof}
By the discussion above, if $a \notin [-1,1]$ then at least one of $b - 2a, b + 2a $ is in $ \C \sm [-2,2]$ for any $b$. Therefore $a \in [-1,1]$ is necessary for $K({ap_n + b})$ connected for $n \geq N$.
But if $a \in (-1,1)$, then $b \pm 2a \in [-2,2]$ for $b$ real and  sufficiently small, more precisely $\vert b \vert \leq 2 ( 1 - \vert a \vert)$. It follows that
\begin{align*}
\lim_{n\rightarrow \infty}\mathcal{C} ( a p_n + b) \subseteq \{ 
(a , b ) \in \mathbb{R}^2  \mid \ 0 <  \vert a \vert \leq 1 ,  \vert b \vert  \leq 2 ( 1 - \vert a \vert) \
 \}.
\end{align*}
For the opposite inclusion, notice that if $(a,b)$ satisfies the right hand side of \eqref{eq:chebyshevconnectedness}, then $ \aphi \circ p_n  ([-2,2]) = \aphi ([-2,2]) \subseteq [-2,2]$ and all critical values $b\pm 2a$ of $\aphi \circ p_n $ are contained in $[-2,2]$. This implies that $J (\aphi \circ p_n  )$ is connected for all $n \geq 2$.
\end{proof}

Thus the connectedness locus for  $a p_n + b$ is asymptotically a diamond shaped region symmetric around the origin in $\mathbb{R}^2$, missing the vertical diagonal.

\subsection{The Orthogonal Sequence of Iterates}

Let $f$ be a monic centered polynomial of degree $d \geq 2$ and let $\mu$ denote the unique measure of maximal entropy. Recall that this measure equals the equilibrium measure for the Julia set $J(f)$ of $f$.

Recall that by Theorem \ref{thm:barnsleyinvariantmeasures} \cite[Thm. 3]{barnsleyinvariantmeasures} the sequence of iterates $\{ f^n \}_n$ is a subsequence of the sequence of orthogonal polynomials.  The following consequence of Corollary \ref{cor:seqofiterates} was suggested by Laura DeMarco.

\begin{theorem}
\label{cor:MMEmain}
 Let $f$ be a monic centered polynomial and let $\{ f^n \}_n$ denote the sequence of iterates. Let $\mu$ be the unique measure of maximal entropy for $f$ and let $\aphi (z) = az +b$.
  Then  
\begin{align}
\label{eq:MMEmain1}
J(f) \subseteq \liminf_{n\rightarrow \infty} J({\aphi\circ f^n}) \subseteq  \limsup_{n\rightarrow \infty} K({\aphi\circ f^n}) \subseteq K(f) .
\end{align}
In particular, if $J(f) = K(f)$, then 
\begin{align}
\label{eq:MMEmain2}
\lim_{n\rightarrow \infty} J({\aphi\circ f^n})  = J(f) .
\end{align}
\end{theorem}

\begin{remark}
Notice that for this special sequence of orthogonal polynomials we obtain a stronger version of Theorem \ref{thm:mainregular}, in particular of Corollary \ref{cor:mainmonic}, since $K(f)$ is convex only if $f$ is a monomial or Chebyshev \cite[Thm. 9.1]{brolin}. 
\end{remark}

\begin{proof}
Since the sequence of iterates is orthogonal by Theorem \ref{thm:barnsleyinvariantmeasures} and the measure of maximal entropy is supported on $J(f)$, it follows directly from Corollary \ref{cor:mainmonic} that
\begin{align*}
J(f) \subseteq \liminf_{n\rightarrow \infty} J({\aphi\circ f^n}) \subseteq  \limsup_{n\rightarrow \infty} K({\aphi\circ f^n} )\subseteq \mathrm{Co}(K(f)) .
\end{align*}

To show the stronger containment we note that if $z_0 \notin K(f)$, then there is an open neighbourhood $U \ni z_0$ such that $\overline{U} \cap K(f) = \emptyset$. Let $\aphi(z) = a z + b$. Choose $R > 0$ such that if $\vv z \vv > R$, then $\vv f(z) \vv > 2 \vv z \vv$ and $N$ so that $\vv a \vv 2^N R + \vv b \vv > R$. 
Then $\vv z \vv > 2^N R$ implies that $\vv \aphi  (z) \vv > R$.  Thus for $\vv z \vv > R$ and $n > N$ we obtain $\vv \aphi\circ f^n (z) \vv > R$ and therefore $K({\aphi\circ f^n (z)}) \subseteq \mathbb{D}(R)$. 

Now choose $N'$ such that $f^{N'} (\overline{U}) \subseteq \C \sm \mathbb{D}(R)$. For $z \in U$ and $n \geq N + N'$ we have that $\vv \aphi \circ f^n (z) \vv > R$ and thus $\overline{U} \cap K({\aphifn}) = \emptyset$. Therefore $z_0 \notin \limsup_{n \rightarrow \infty} K({\aphi \circ f^n})$ and hence
\begin{equation*}
\limsup_{n \rightarrow \infty} K({\aphi \circ f^n}) \subseteq K(f),
\end{equation*}
yielding \eqref{eq:MMEmain1}.  Note that \eqref{eq:MMEmain2} follows from \eqref{eq:MMEmain1} by definition of $\liminf$ and $\limsup$. 
\end{proof}

\subsubsection{Connectedness Loci}

Suppose  $f$ is a monic centered  polynomial where all critical points of $f$ belong to the filled Julia set $K(f)$, so that $K(f)$ is connected.  And let $\aphi(z) = a z +b$.
Then the critical points and critical values of $f^n$ are in $K(f)$ by total invariance.
 The critical points of $\aphifn$ coincides with the critical points of $f^n$. By the chain rule $c_i$ is a critical point of $f^n$ if and only if $f^k (c_i)$ is a critical point of $f$ for some $k$ with $0 \leq k < n$.
Similarly, $z$ is a critical value of $f^n$ if and only if $ z = f^k (c_i)$ for some critical point $c_i$ and $1 \leq k \leq n$.
 If  $z$ is a critical value of $f^n$, then the corresponding critical value of $\aphifn$ is $\aphi (z)$. 
In this section we will explore the connectedness loci of $\aphifn$ as $n$ tends to infinity. 

 In light of  Theorem \ref{cor:MMEmain} we  have the following upper bound on the connectedness locus.

\begin{proposition}
Let $c$ be a critical point of $f$. If there is $k$ such that $\aphi \circ f^{k} (c)\notin K(f)$, then there is $N$ so that $\aphi \circ f^{k} (c) \notin K({\aphi \circ f^{n}})$
 for $n \geq N $. Thus $(a,b) \notin \displaystyle\limsup_{n \rightarrow \infty} \mathcal{C} (\aphifn) $.
 \end{proposition}
 \begin{proof}
 $\aphi \circ f^{k} (c)\notin K(f)$ implies that $\aphi \circ f^{k} (c)\notin \limsup_{n \rightarrow \infty} K({\aphifn})$ by \eqref{eq:MMEmain1}, so there is $N$ such that $\aphi \circ f^{k} (c)\notin  K({\aphifn})$ for $n \geq N$. 
That is, at least one critical value eventually escapes 
and thus the Julia sets are disconnected for sufficiently large $n$. 
 \end{proof}

Let $C(f)$ denote the set of critical points of $f$ and $V(f)$ the set of critical values.

Suppose $z_i \in V(f) \cap J(f)$. The orbit of this critical value is very sensitive to small perturbations by $\aphi$. 
  Therefore, here we will  consider only the case  where all critical points and critical values of $f$ are contained in the interior of $K(f)$. We will moreover restrict to hyperbolic $f$ since in the parabolic case a   point in a parabolic basin may shift to the basin of infinity  under arbitrary small perturbations.

So suppose in the following that $f$ is a monic centered hyperbolic polynomial with all critical points in the interior of the filled Julia set. 
Note that in the following $d(X,Y) = \inf_{x \in X, y\in Y} \vv x - y \vv$, while 
 $d_h$ denotes the asymmetric  Hausdorff semi-distance $d_h(X, Y) = \sup_{x\in X} \inf_{y \in Y} \vv x - y \vv $.

\begin{lemma} If $f$ is a monic centered  hyperbolic polynomial, then
there is a positive integer $N$ and an open set $U$ in the parameter space around $(1,0)$   so that $U \subseteq  \mathcal{C} (\aphifn)$ for $n \geq N$, i.e. $U \subseteq \liminf_{n \rightarrow \infty}  \mathcal{C} (\aphifn)$.
\end{lemma}

\begin{proof}
Let $Z = \{ z \mid \exists k \text{ such that } f^k (z) = z , \vv( f^k ) ' (z)  \vv < 1 \}$ denote the set of attracting periodic points. 
Choose for each periodic attracting point $z\in Z$  a disc $D_z$   
centered at $z$  such that $\overline{f(D_z)} \subseteq D_{f(z)}$, i.e. ${f(D_z)}$ is compactly contained in $D_{f(z)}$.
Define $D_n :=  \bigcup_{z \in Z} f^{-n} (D_z) = f^{-n} ( \bigcup_{z \in Z} D_z )$. Then $\overline{f (D_n)} = \overline{D_{n-1}} \subseteq D_n$.
As $f$ is hyperbolic by assumption all critical points are attracted to some periodic attracting point $z \in Z$, so there is $N$ such that $C(f) \subseteq D_N$. 
It follows that $V(f^n) = \bigcup_{k=1} ^n f^k (C) \subseteq 
D_{N-1}$ for $n \geq N$. 

Choose $R$ such that $K(f) \subseteq \mathbb{D}(R)$.
If $\aphi$ is sufficiently close to the identity in $\mathbb{D}(R)$,
 i.e. if  $\vv\vv \aphi (z) - z \vv\vv_{\mathbb{D}(R)} < \epsi$ for $\epsi = d (D_{n+1} , \partial D_n )$, then $V(\aphifn) \subseteq 
\aphi( D_{N-1} )$ for $n \geq N$. Hence the orbit of any critical point of $\aphifn$ is bounded  and thus contained in $K({\aphifn})$. 
That is, all $(a,b)$ such that $\vv\vv \aphi (z) - z \vv\vv_{\mathbb{D}(R)} < \epsi$ are contained in $ \mathcal{C} (\aphifn)$ and hence there is an open set around $(1,0)$ contained in the connectedness locus. 
\end{proof}

A natural question is then how big this open set might be. The following propositions explore this.

\begin{proposition} If $f$ is a monic centered hyperbolic polynomial and
if there is a compact subset $L \subseteq \overset{\circ} {K}_f$ such that $\bigcup_{n \geq 1} \bigcup_{c \in C} \{ \aphifn (c) \} \subseteq L$, then there is an  $N ' \in \mathbb{N}$ such that $(a,b) \in \mathcal{C} (\aphifm)$ for all $m \geq N '$, i.e. $(a,b) \in \liminf_{m \rightarrow \infty} \mathcal{C} (\aphifm)$. 
\end{proposition}

\begin{proof}
As above, let $D_n = f^{-n} (\bigcup_{z \in Z} D_z)$ for $n \in \mathbb{Z}$ and $D_z$ the discs centered at the attracting periodic points in the interior $\overset{\circ}{K}_f$. Then $\overset{\circ}{K}_f = \bigcup_{n\geq 0} D_n$.  Suppose $L$ exists as prescribed. As  $L$ is compact there is $N$ such that $ L \subseteq \overline{D_ N} \subseteq D_{N+1}$. 
Thus if $\aphifm (c) \in L$ for all $m\geq 1$ and all $c \in C(f)$, then 
$\aphifm (c) \in D_{N+1}$ for all $m \geq 1$ and all $c \in C(f)$.

$L$ is compact, so  $ \overline{ \bigcup_{n \geq 1} \bigcup_{c \in C} \{ \aphifn (c) \} } \subseteq L$. 
As any periodic attracting point is a limit point of a critical orbit 
 it follows that $\aphi(Z) \subseteq L$. 
Consider $\aphifm (D_{N+1}) = \aphi(D_{N+1 - m}) $. 
Let $\epsi$ be the distance $d ( D_N , \partial D_{N+1})$. 
There is  a natural number $M$ such that $d _H(D_{-m}, Z) < \frac{\epsi}{\vv a \vv}$ for $m > M$. If $m > M + N + 1$, then $D_{N+1 - m} \subseteq D_{-M}$, thus $\aphifm (D_{N+1}) \subseteq \aphi (D_{-M})$. 
Since $d_H (\aphi (D_{-M}) ,\aphi(Z) ) < \epsi $ by the choice of $M$, it follows that $\aphifm (D_{N+1}) \subseteq  D_{N+1} $  for $m > N + M + 1 = N'$. So evidently $D_{N+1}$ is a trapping region containing all critical values, so no critical values escape for $m > N'$.
\end{proof}

The critical values of the tail of the sequence $\aphifn$ are controlled as the critical points of $f$ are attracted to attracting cycles. This leads to the following.

\begin{proposition} Suppose $f$ is a monic centered hyperbolic polynomial and
 $\aphi (Z) \in \overline{{D}_N}$ for some $N$, where $Z$ is the set of  attracting periodic points. Then there exists $N' = N' (N)$ such that if  $$\bigcup_{1 \leq n \leq N'} \bigcup_{c \in C(f)} \{ \aphifn (c) \} \subseteq \overline{{D}_N},$$ then  $(a,b) \in \mathcal{C}(\aphifk)$ {for  $k$ sufficiently large.}
\end{proposition}
The quantification on $k$ is vague, but we note that the condition is on a finite number of iterates.
\begin{proof}
Assume $\aphi (Z) \subseteq \overline{D_N} \subseteq D_{N+1}$ for some $N$. 
As above, let $\epsi = d (D_N , \partial D_{N+1} )$. 
There is $M$ such that that $d_H (\overline{D_{-M}}, Z) < \frac{\epsi}{\vert a \vert}$.  
There is $N'$ so $f^k (c) \in D_{-M}$ for any $c \in C(f)$ and $k > N'$.  In particular,  $\aphifk (c) \in \aphi(D_{-M}) \subseteq D_{N+1}$ for $k > N'$. Recall that $D_{N+1}$  is a trapping region containing all critical values for $k$ sufficiently big
 by the proof above. 
Thus if   $\bigcup_{1 \leq n \leq N'} \bigcup_{c \in C(f)} \{ \aphifn (c) \} \subseteq \overline{{D}_N} \subseteq D_{N+1}$,
then all critical values are contained in the trapping region $D_{N+1}$, 
so  $(a,b) \in \mathcal{C}(\aphifk)$ for $k$ sufficiently big. 
\end{proof}

\section*{Acknowledgements}

The second author would like to thank the Danish Council for Independent Research \textbar\  Natural Sciences for support via the grant DFF – 4181-00502. The second author would like to thank Laura DeMarco for helpful conversations. This paper is inspired by the master’s thesis written by the first author at IMFUFA, Institute of Science and Environment at Roskilde University. The authors would like to thank IMFUFA for its hospitality during the conception and initial writing of the paper.

\end{document}